\documentclass{article}

\usepackage{arxiv}

\usepackage[utf8]{inputenc} 
\usepackage[T1]{fontenc}    
\usepackage{hyperref}       
\usepackage{url}            
\usepackage{mathrsfs}
\usepackage{amsmath}
\usepackage{graphicx}
\usepackage{cite}
\usepackage{caption}
\usepackage{subcaption}
\usepackage{amssymb}
\usepackage{algorithm}
\usepackage[noend]{algpseudocode}
\usepackage[short]{optidef}

\title{On the benefit of overparameterisation in state reconstruction: \\ An empirical study of the nonlinear case}

\author{Jonas F.\ Haderlein\thanks{This research was conducted in the Australian Research Council Training Centre in Cognitive Computing for Medical Technologies (project number ICI70200030) and funded by the Australian Government.},~Andre D.\ H.\ Peterson%
\thanks{A.D.H.\ Peterson and D.B.\ Grayden are with the Department of Medicine, St Vincent's Hospital, The University of Melbourne.},~Parvin Zarei Eskikand,
\AND Anthony N.\ Burkitt,~Iven M.\ Y.\ Mareels\thanks{I.\ Mareels is with the Institute for Innovation, Science and Sustainability, Federation University Australia, Mt Helen, Vic 3350},~David B.\ Grayden
\thanks{J.\ Haderlein, A.D.H.\ Peterson, P.\ Zarei Eskikand, A.N.\ Burkitt and D.B.\ Grayden are with the Department of Biomedical Engineering, University of Melbourne VIC 3010, Australia.}
\thanks{D.B.\ Grayden is with the Graeme Clark Institute, University of Melbourne VIC 3010, Australia.}
}

\begin{document}

\maketitle

\begin{abstract}

The empirical success of machine learning models with many more parameters than measurements has generated an interest in the theory of overparameterisation, i.e., underdetermined models.
This paradigm has recently been studied in domains such as deep learning, where one is interested in good (local) minima of complex, nonlinear loss functions.
Optimisers, like gradient descent, perform well and consistently reach good solutions.
Similarly, nonlinear optimisation problems are encountered in the field of system identification. Examples of such high-dimensional problems are optimisation tasks ensuing from the reconstruction of model states and parameters of an assumed known dynamical system from observed time series.
In this work, we identify explicit parallels in the benefits of overparameterisation between what has been analysed in the deep learning context and system identification.
We test multiple chaotic time series models, analysing the optimisation process for unknown model states and parameters in batch mode. We find that gradient descent reaches \textit{better} solutions if we assume more parameters to be 
\textit{unknown}.
We hypothesise that, indeed, overparameterisation leads us towards better minima, and that more degrees of freedom in the optimisation are beneficial so long as the system is, in principle, observable.

\end{abstract}


\section{Introduction}

The theoretical analysis of deep learning has recently observed unexpected behaviour in the training of very large, overparameterised models \cite{Belkin2019, Liu2020a}: modern state-of-the-art models are often underdetermined (i.e., they contain more unknown parameters than available measurements) but perform empirically well. This stands in stark contrast to classical beliefs concerning the bias-variance trade-off: models that interpolate the measurements (small bias) were assumed to perform poorly on out-of-sample data (high variance).
Spawning from observations that apparently contradicted this paradigm, a theory of overparameterisation has started to develop \cite{Belkin2019, Gerace2020a, Hastie2019}. 
More research into the benefits of overparameterisation has since been conducted, primarily in the deep learning context \cite{Liu2021, Liu2022}.

The training of (deep) neural networks is a nonlinear optimisation problem.
It is usually performed by some version of gradient descent, optimising the weights with respect to a loss function (such as the mean-squared error) and the target output of the system (usually given by the measurements). 
A key to the development of a theory for overparameterised models was the analysis of the ensuing loss landscape, which has been found to exhibit properties that make the search for (global) optima with gradient descent very efficient \cite{Liu2020a, Liu2022}. 
In particular, it has been found that these problems are not locally convex, but instead satisfy (almost everywhere) a version of the Polyak-Lojasiewicz condition, which makes it easy to reach good minima via gradient descent. 
Similar tendencies have also been noted in the training of other models in the machine learning community \cite{Xu2018, Hastie2019, Buhai2020}, independent of the training algorithms.

Complex optimisation tasks also arise naturally in the field of system identification. Here, one is interested in the reconstruction of a system state from measurements and a postulated model family parameterised by model coefficients.
The reconstruction of such unknown model states and parameters is a problem approached with a multitude of different approaches \cite{Schoukens2019, Bittanti2019}. Common to these approaches is the minimisation of a loss function, usually defined by both measurements and the prior assumption on the model class (usually in the form of a dynamical system). 
Under some generic assumptions of observability, a batch optimisation of high-dimensional time-varying model states and parameters has been shown to work well \cite{Mareels2002, Haderlein1}. 
This is, in some sense, surprising since measurements are typically lower dimensional than the full list of time-dependent states and parameters characterising the model. 

It seems that state reconstruction has some unexpected parallels to the aforementioned overparameterised training of machine learning models. In this work, we analyse, empirically, whether there are, in fact, benefits of overparameterisation in state reconstruction. 
In order to provide complex but well-posed problems, we analyse a class of nonlinear dynamical systems exhibiting chaotic behaviour. 
We solve for a vector of unknowns that is larger, by design, than the vector of measurements. 
How does such an optimisation compare to a hypothetical scenario in which we have fewer unknowns due to prior knowledge? Does gradient descent find good state and parameter reconstructions or does it get stuck in local, spurious minima? We first review relevant literature and then provide an illustrative example.

\paragraph{State reconstruction as optimisation problem}

The reconstruction of states and parameters of autonomous dynamical systems from noisy measurements is, in general, a non-trivial task. There is a vast literature dealing with this fundamental problem, briefly summarised in the following.

In time series state reconstruction and denoising, methods usually focus on finding an initial state from which the whole trajectory can be recreated (such as \cite{Ortega2015a}). 
The complexity of the ensuing loss function in such a simulation error framework has, for example, been reported by \cite{Voss2004}. 
Methods exist to find solutions through recursions \cite{Bittanti2019, Anderson1990}.

In case we observe the full state trajectory and aim to recover parameters only (as is, for example, assumed in autoregressive modelling), one can often rephrase the optimisation as a regression problem, for which many known techniques exist \cite{Schoukens2019, Brunton2016b}. 
This approach, however, suffers from an inevitable curse called `errors-in-variables' when the measurements are corrupted by noise \cite{Voss2004}.

Methods have been proposed to linearise certain system identification tasks via the introduction of additional regression parameters, a type of overparameterisation that has been found to have negative effects \cite{Ortega2020, Cisneros2022a}. 
We argue that this parameterisation is qualitatively different to the one in this work, where we restrict ourselves to cases where the structure of the model is known.
We perform a direct but nonlinear reconstruction of both unknown states and parameters that are actually present in the true model. 
We show that gradient descent robustly converges towards solutions close to the ground truth in this case -- counter-intuitively so.

\paragraph{Overparameterisation in optimisation theory}

The phenomenon of overparameterisation has been extensively studied from the viewpoint of machine learning. 
Here, the classical bias-variance trade-off has been subsumed by a so-called `double descent' \cite{Liu2020a, Hastie2019, Liu2021, Ribeiro2021, Gerace2020a}, describing the model's ability to generalise to new data.
A second `descent' of risk appears beyond the interpolation point, i.e., in the overparameterised regime \cite{Belkin2018, Belkin2019}.

Gradient descent-based optimisers are usually employed for finding solutions to high-dimensional learning tasks.
The study of such systems is thus closely coupled to a theory of gradient descent in high dimensions \cite{Saxe2014, Gerbelot2022, Advani2020} and complex nonlinear optimisation tasks \cite{Mannelli2020b}, observing a property of high-dimensional gradient flow to find good non-spurious minima. 
We empirically explore high-dimensional gradient flow for loss functions structurally somewhat different to those in machine learning, yet find similar tendencies. 

\section{A motivating example}
\label{henonexample}
To illustrate the complex loss functions potentially arising from the optimisation of model states, we demonstrate an example from the H\'enon map, here represented 
 as a second order autoregressive difference equation:

\begin{equation}
\begin{array}{lcll}
   x_{t+1} &=& {\Theta}_1+  {\Theta}_2 x_{t}^2 + {\Theta}_{3} x_{t-1}, \\
   y_t &=& x_t + \mu_t,
   \end{array}  \label{henon}
\end{equation}

with (observed) variable $x_{t}$ and known model parameters $\Theta_1 = 1, \Theta_2 = -1.4, \Theta_{3} = 0.2$ and the discrete time index $t = 1, \cdots, N$. Given a random $x_{1}, x_{2}$, we simulate $x_t$ trajectories from equation (\ref{henon}). 
Note that the parameters are chosen so that the iteration is chaotic, i.e., nearby initial conditions diverge quickly.

In this example, we assume the measurement noise $\mu_t = 0 ~ \forall t$; we denote $X_N=(x_2, \cdots x_N)^T = Y_N = (y_2, \cdots y_N)^T$ the states except for the initial $x_1$. 
In this artificial case, we show our inability to reconstruct the time series $Y_N$ as a function of the initial value only: Finding the original $x_1 = 0.629345\ldots$ from knowledge of $x_2 = 0.450339\ldots$ to $x_N$ by simulation starting with $\hat{x}_{1} = x_{1} + \epsilon $ for $\epsilon = -0.01, \cdots, 0.01$, resulting in $\hat{X}_N = (\hat{x}_2, \cdots \hat{x}_N)^T$. 

\begin{figure}[]
     \centering
     \begin{subfigure}{0.4\columnwidth}
    \includegraphics[width=\columnwidth]{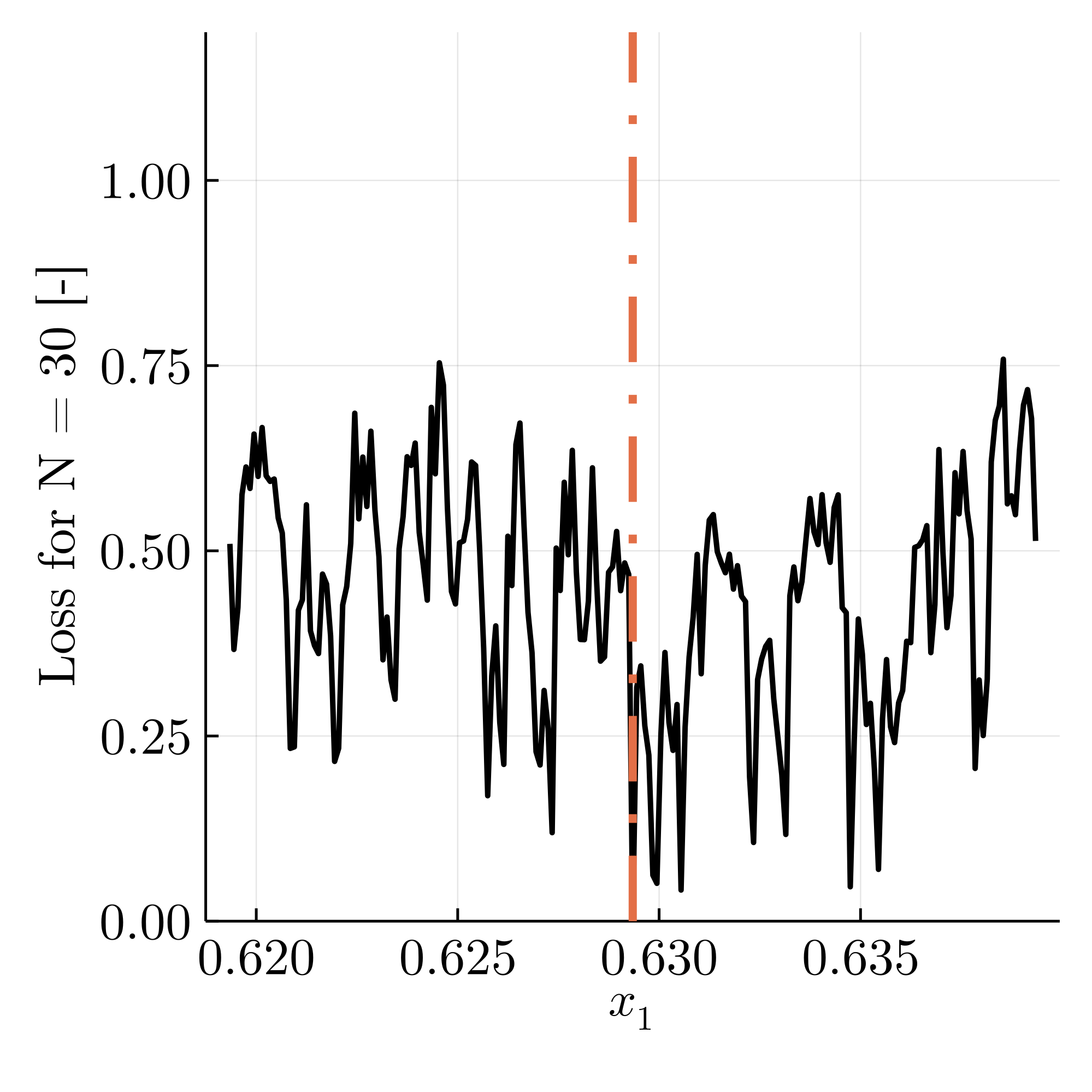}
         \caption{}
     \end{subfigure}
     \hspace{1mm}
     \begin{subfigure}{0.4\columnwidth}
    \includegraphics[width=\columnwidth]{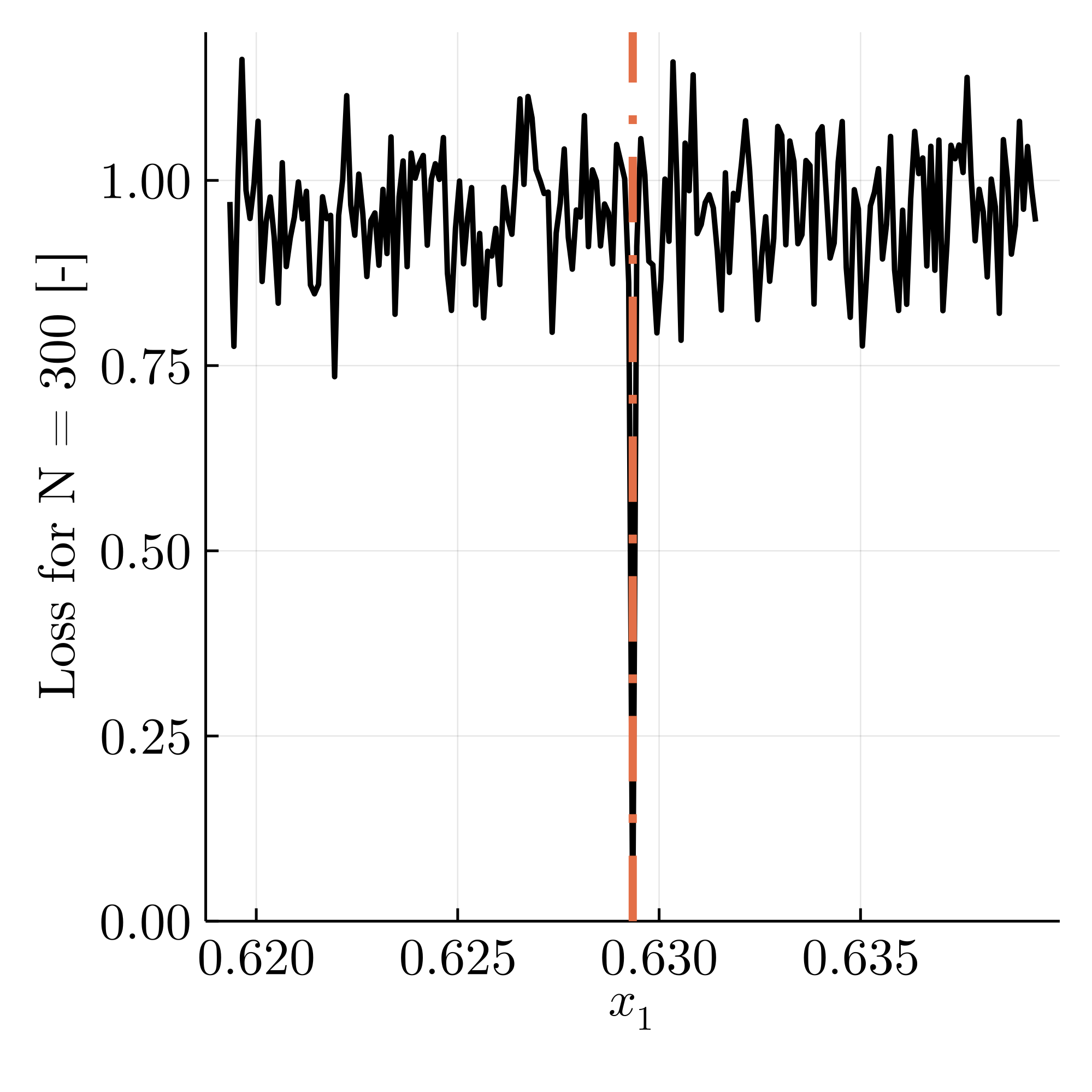}
         \caption{}
     \end{subfigure}
    \caption{Loss $\frac{1}{N} \ell_0$ over the initial value of $\hat{x}_{1}$ for (a) $N = 30$ and (b) $N = 300$. The randomly generated true initial condition ${x}_{1} = 0.629\ldots$ is marked with an orange line.}
    \label{figure1}
\end{figure}

The aim is to minimise the error $\ell_0 = \|X_N-\hat{X}_N\|^2_2  =  \sum_{t=2}^{N}\| y_t - \hat{x}_t \|_2^2$.
This procedure is equivalent to the optimisation program
\begin{mini}
  {\hat{x}_{1}}{\frac{1}{N} \ell_0(Y_N, \hat{x}_1)~~~~~~~~~~~~~~~~~~~~~~~~~~~~~~~}{}{}
  \addConstraint{\hat{x}_{t+1}}{= {\Theta}_1+  {\Theta}_2 \hat{x}_{t}^2 + {\Theta}_{3} \hat{x}_{t-1}, ~ \hat{x}_2 = x_2}.
 \label{simulationerror} 
 \end{mini}
 We plot the loss for $X_{30}$ and $X_{300}$, respectively, over the initial deviation; see Figure \ref{figure1}. 

Simulating trajectories forward in time and optimising only for an initial condition is called the `simulation-error' approach \cite{Schoukens2019}. It becomes clear that, for large enough $N$, such an optimisation of $\hat{x}_{1}$ alone is likely to only reach a (spurious) local minimum in a search for the true initial value $x_{1}$. In case the system is perturbed by a small noise, the shadowing theorem \cite{Ott2002} indicates that the noisy trajectory is actually close to a true, not noisy, trajectory, hence making the exact recovery of $x_1$ impossible using the above loss function.

In the following, we consider multiple relaxations of the extreme case in the program (\ref{simulationerror}).

\section{Methods}

\subsection{Loss functions}

Consider a 
state space map over a finite time horizon $t = 1, \cdots, N$,
\begin{equation}
\label{statefunction}
\begin{array}{lcll}
   x_{t+1} &=& f(x_t, \Theta) + \nu_t, & x_1 = \kappa, \\ 
   y_t &=& h(x_t,\Theta) + \mu_t,
   \end{array}
\end{equation}
where $f$ is the state transition function, $x_t \in R^n$ is the state at time $t$, with parameters $\Theta \in R^m$, and $h$ is the measurement function, measurements $y_t \in R^p$, measurement noise $\mu_t$, and 
state transition noise $\nu_t$. 

Denote $X=(x^T_1, \cdots x^T_N)^T$ and $Y=(y^T_1, \cdots y^T_N)^T$.

\paragraph{Relaxation 1} First, the aim is to reconstruct a denoised $X$ from knowledge of $Y$ as well as $f, h$ and $\Theta$. 
To this end, we reconstruct all state variables by adding the terms for each state transition over $t$ (i.e., incorporating the uncertainty $\nu_t$):
\begin{mini}
  {\hat{\mathcal{X}}_1}{\frac{1}{N}
  \ell_1(Y, \Theta, \hat{X})
  }{}{},
 \label{relaxation1} 
 \end{mini}
 with 
 \begin{equation}
 \begin{aligned}
     \ell_1(Y, \Theta, \hat{X}) = &\sum_{t=1}^{N-1} \| \hat{x}_{t+1}-f(\hat{x}_t,\Theta)\|_2^2 + \\ &\rho \sum_{t=1}^{N}\| y_t - h(\hat{x}_t, \Theta) \|_2^2,
\end{aligned}
 \end{equation}
with estimates $\hat{\mathcal{X}}_1 = \hat{{X}} \in R^{N \times n}$ and where $\rho > 0$ is a scalar controlling our prior belief of the relative magnitudes of measurement and state transition error.

\paragraph{Relaxation 2}

Second, if $\Theta$ is unknown, we additionally optimise for $m$ parameter estimates $\hat{\Theta}$ , so that
\begin{mini}
  {\hat{\mathcal{X}}_2}{\frac{1}{N} \ell_2(Y, \hat{\Theta}, \hat{X})}{}{},
 \label{relaxation2} 
 \end{mini}
  with 
 \begin{equation}
 \begin{aligned}
     \ell_2(Y, \hat{\Theta}, \hat{X}) = &\sum_{t=1}^{N-1} \| \hat{x}_{t+1}-f(\hat{x}_t,\hat{\Theta})\|_2^2 + \\ &\rho \sum_{t=1}^{N}\| y_t - h(\hat{x}_t, \hat{\Theta}) \|_2^2,
\end{aligned}
 \end{equation}
where $\hat{\mathcal{X}}_2 = (\hat{{X}}^T, \hat{\Theta}^T)^T \in R^{(N \times n) + m}$.

This case has some significance since, in practice, it is the reconstruction of unknown system parameters that we are most interested in.

\paragraph{Relaxation 3}
Third, we consider the case of potentially time-varying parameters, so that the state map reads 
\begin{equation}
\label{augmentedstatefunction}
\begin{array}{lcll}
   x_{t+1} &=& f(x_t, \theta_t) + \nu_t, & x_1 = \kappa, \\ 
   \theta_{t+1} &=& \theta_t + \vartheta_t , & \theta_1 = \Theta, \\
   y_t &=& h(x_t,\theta_t) + \mu_t,
   \end{array}
\end{equation}
where $\theta = (\theta^T_1, \cdots, \theta^T_N)^T \in R^{N \times m}$ are time-varying, unknown parameters perturbed by zero-mean $\vartheta_t$. 
With parameters effectively being slowly time-varying states, we denote the augmented state vector ${\mathcal{X}}_3 = (x_1^T, \theta_1^T, \cdots, x_N^T, \theta_N^T)^T \in {R}^{N \times (n+m) }$. The program to reconstruct an estimate $\hat{\mathcal{X}}_3$ is 
\begin{mini}
  {\hat{\mathcal{X}}_3}{ \frac{1}{N} \ell_3(Y, \hat{\mathcal{X}_3})}{}{},
 \label{relaxation3} 
 \end{mini}
   with  \cite{Mareels2002, Haderlein1}
 \begin{equation}
 \begin{aligned}
     \ell_3(Y, \hat{\mathcal{X}_3}) = &\sum_{t=1}^{N-1} \| \hat{x}_{t+1}-f(\hat{x}_t,\hat{\theta}_t)\|_2^2 + \\&\rho \sum_{t=1}^{N}\|y_t- h(\hat{x}_{t}, \hat{\theta}_{t})\|^2_2 
    + \delta \sum_{t=1}^{N-1} \|\hat{\theta}_{t+1} -\hat{\theta}_{t}\|^2_2,
\end{aligned}
 \end{equation}
where $\delta > 0$ is a scalar controlling the stationarity of parameters.

\paragraph{Comparison}

Every relaxation brings in more parameters to the optimisation program, thus transferring the loss function into a higher dimensional space. 
Vice versa, each relaxation collapses to its original form without the additional degrees of freedom.
In particular, Relaxation 2 and 3 can be transformed into relaxation 1 in the case where the true parameters are known. 
Hence, Relaxation 2 and 3 are denoising and parameter recovery tasks, whereas Relaxation 1 is purely denoising. 
Relaxation 2 (i.e., optimising for a single constant vector $\hat{\Theta}$ instead of time-varying $\hat{\theta}_t$) is classically used in systems identification \cite{Schoukens2019}. 

We study the loss landscapes of $\ell_{1,2,3}$, comparing the convergence towards the ground truth and potential local minima.
We note that the existence of a unique global minimum as well as the shape of the loss landscape for either relaxation clearly depend heavily on the choice of $f, h$. We show results for three examples of chaotic time series, assuming observability of $X$ or $X,\Theta$ from $Y$ throughout.

\subsection{Optimisation}

We minimise $\ell$ via a (primitive) gradient descent implementation
\begin{equation}
\label{GD}
    \hat{\mathcal{X}}^{k+1} \leftarrow \hat{\mathcal{X}}^{k} - \eta \frac{\nabla \ell(\hat{\mathcal{X}}^k)}{\|\nabla \ell(\hat{\mathcal{X}}^k) \|_2}, ~~ k = 1, \cdots, K, 
\end{equation}
where $k$ is an index denoting the iteration of gradient descent steps, $\nabla \ell(\hat{\mathcal{X}})$ denotes the derivative of $\ell$ at $\hat{\mathcal{X}}$, and $\eta$ is a learning rate.
Unless otherwise stated, we start the optimisation from a random initial value for $\hat{\mathcal{X}}^1 \sim \mathcal{N}(0,1)$\footnote{This notation means that every component in the vector $\mathcal{X}$ is normally distributed with zero mean and variance 1.}and use a small learning rate of $\eta = 0.01$\footnote{Scaling is important in this context, consider $\eta$ small relative to the variance in the random initial value.}  to give a proper indication of the shapes of the loss landscapes. 

For each $k$, we monitor the loss $\ell_{1,2,3}$ as well as the state reconstruction error $e(\hat{\mathcal{X}}) = \|\hat{\mathcal{X}}-\mathcal{X}\|_2 / \|\mathcal{X}\|_2$ with ground truth ${\mathcal{X}} = (x_1^T, \theta_1^T , \cdots, x_N^T, \theta_N^T)^T$ and $\hat{\mathcal{X}}$ being the different outcomes of Relaxations 1, 2, or 3 in the same format -- thus, in case 1, all $(\hat{\theta}_1^T , \cdots, \hat{\theta}_N^T)^T$ are filled in with the ground truth, and in case 2 with the single estimate $\hat{\Theta}$.

The task of performing the above algorithm for large $N,n,m$ is, in practice, possible via modern automatic differentiation frameworks like \textit{Zygote} \cite{Innes2018} in Julia, which has been used in this work.

\section{Numerical Examples}
\label{Examples}
In this section, we compare the fully overparameterised loss with all states and parameters assumed unknown (Relaxation 3) versus the seemingly simpler Relaxations 1 and 2 in which some unknowns are replaced by the exact, ground truth values. 
We simulate in total ten (an arbitrary number) trajectories of each chaotic map starting from a random initial condition $\kappa \sim \mathcal{N}(0,1)$, unless stated otherwise, for which the simulated time series does not diverge nor collapse. 

We use a time horizon of $N = 300$ and select $\delta = 10$, and $\rho = 0.1$ or $1$, as indicated in the examples.

We indicate the complexity for each example, including the scaling with $N$, in Table \ref{TableProj2bRes} at the end of the section.

\subsection{State recovery of the H\'enon map}

The task here is to find a state and parameter estimate for a H\'enon map time series simulated from equations (\ref{henon}) where $\mu \sim \mathcal{N}(0,\sigma)$ is additive measurement noise. The loss function for Relaxation 3 in this case reads 
\begin{equation}
\label{henonloss}
\begin{aligned}
    \ell_3^{H}(\hat{\mathcal{X}}_3) =&\sum_{t=2}^{N-1} \| \hat{x}_{t+1}-(\hat{\theta}_{1,t} + \hat{\theta}_{2,t} \hat{x}^2_{t} + \hat{\theta}_{3,t} \hat{x}_{t-1})\|_2^2 + \\ &\rho \sum_{t=1}^{N}\|y_t-\hat{x}_{t}\|^2_2 
    +\delta \sum_{t=1}^{N-1} \|\hat{\theta}_{t+1} -\hat{\theta}_{t}\|^2_2.
\end{aligned}
\end{equation}
We infer knowledge about $n+m = 4$ states, while provided with $p = 1$ measurements for each time step $t$. 

Loss $\ell_2^{H}(\hat{\mathcal{X}}_2)$ can directly be derived from the equation (\ref{henonloss}) by replacing all time-varying $\hat{\theta}_t$ by a corresponding single estimate $\hat{\Theta}$.  
Relaxation 1, i.e., assuming $\Theta_1,\Theta_2, \Theta_{3}$ to be known, can be derived from the equation (\ref{henonloss}) by replacing $\hat{\theta}_t$ by their corresponding ground truth values. 

\begin{figure}[]
     \centering
     \begin{subfigure}{0.35\columnwidth}
    \includegraphics[width=\columnwidth]{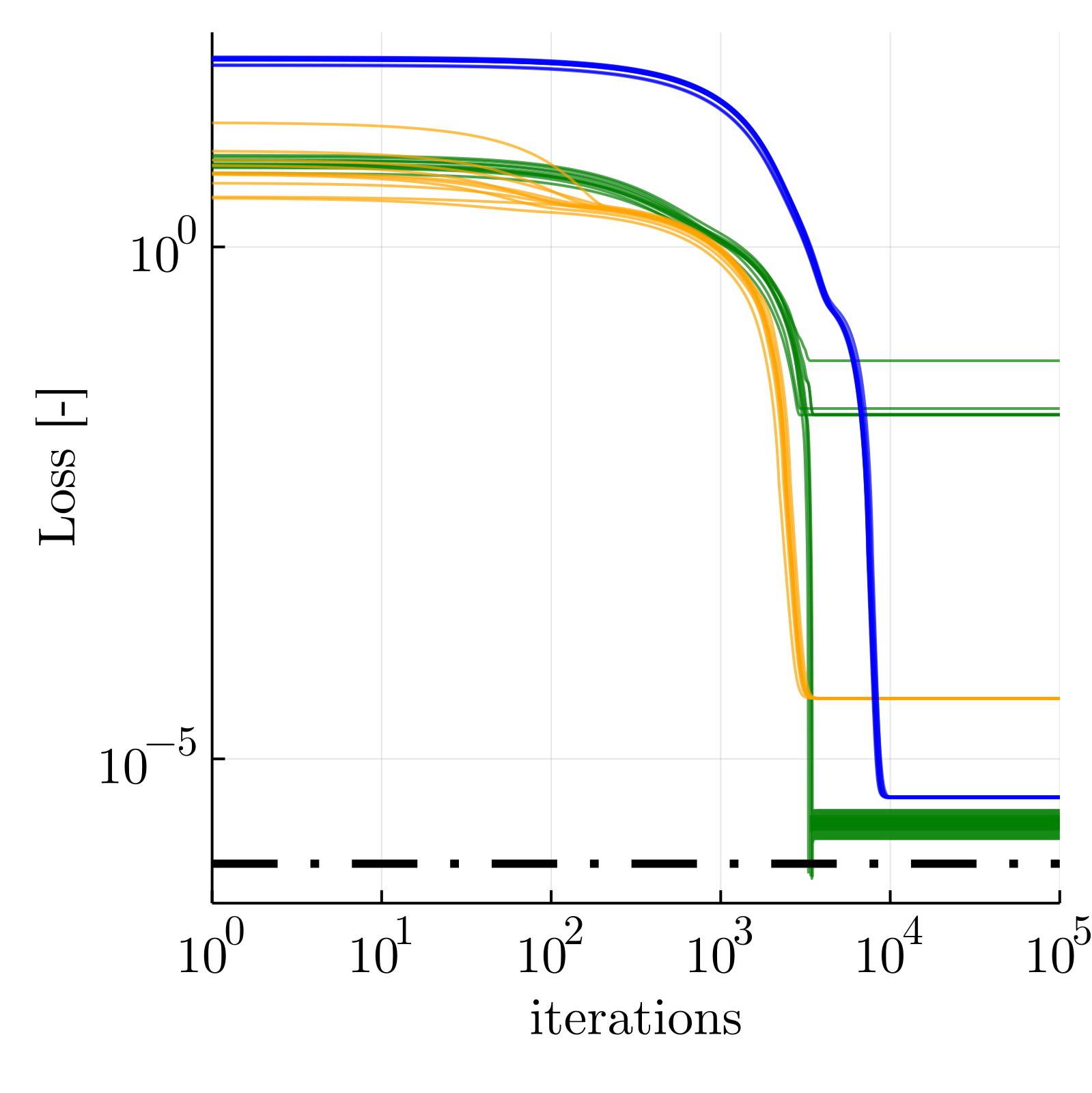}
         \caption{}
     \end{subfigure}
     \hspace{1mm}
     \begin{subfigure}{0.35\columnwidth}
    \includegraphics[width=\columnwidth]{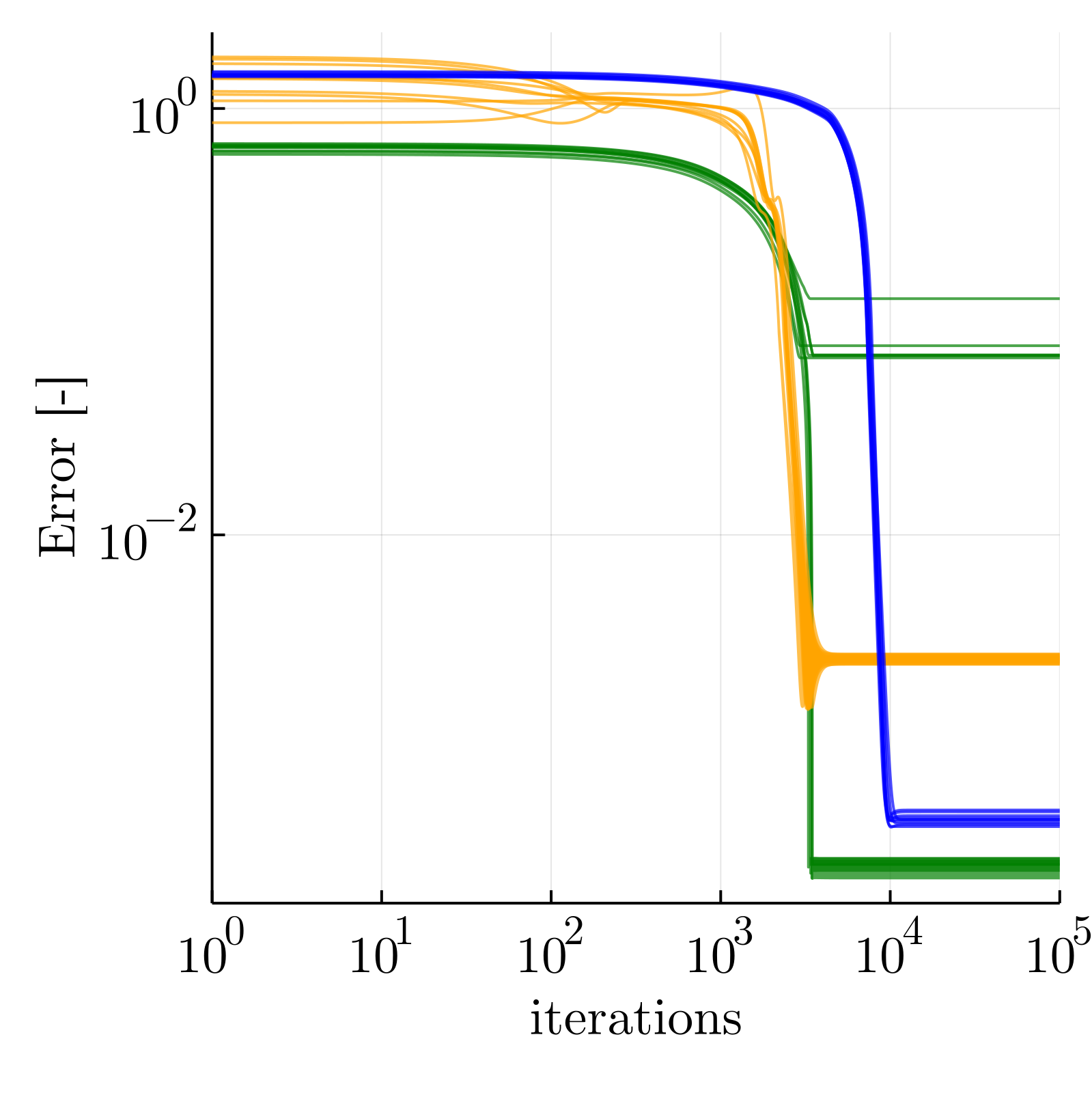}
         \caption{}
     \end{subfigure}
\vspace{1mm}
     \begin{subfigure}{0.35\columnwidth}
    \includegraphics[width=\columnwidth]{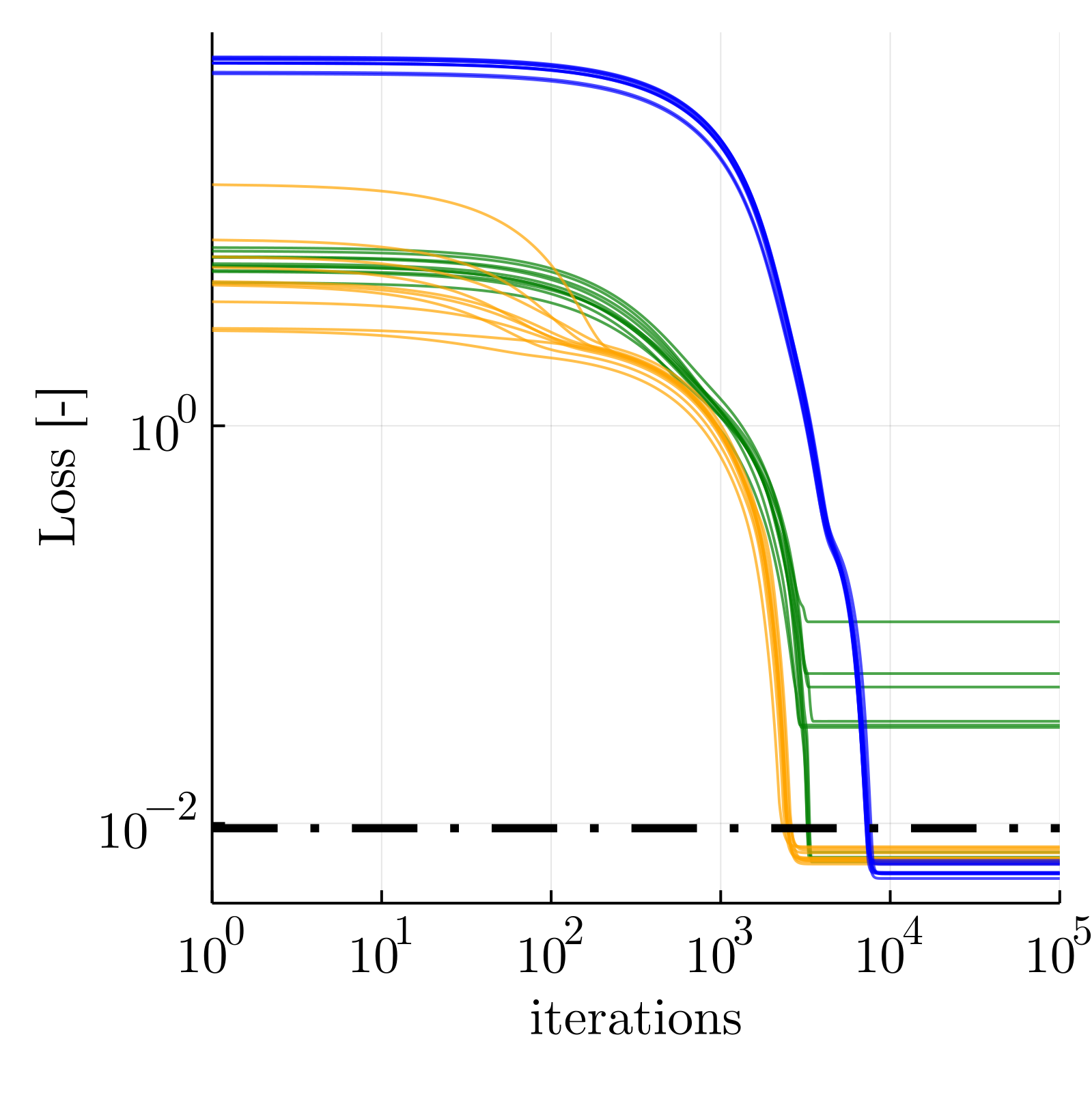}
         \caption{}
     \end{subfigure}
     \hspace{1mm}
     \begin{subfigure}{0.35\columnwidth}
    \includegraphics[width=\columnwidth]{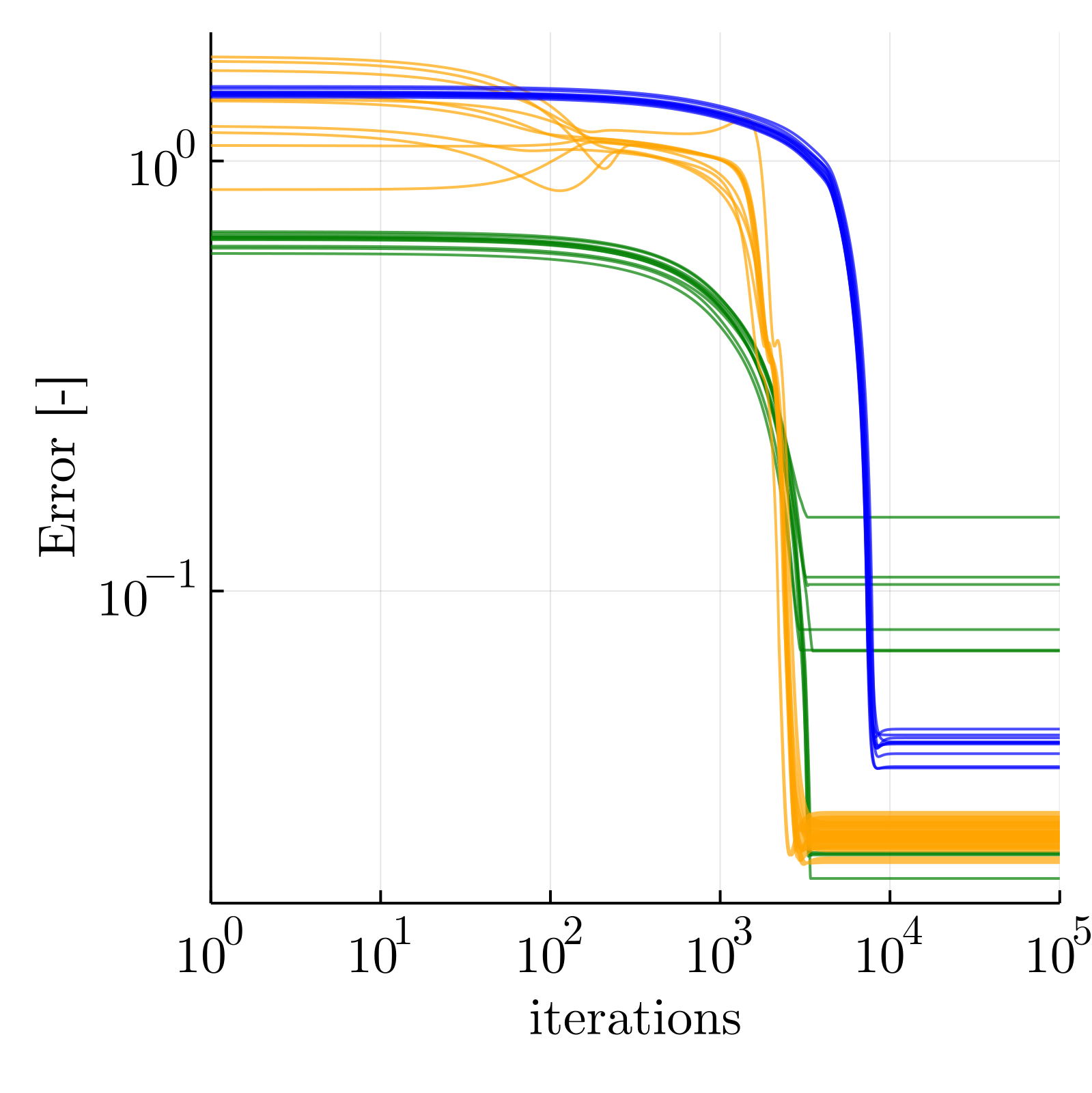}
         \caption{}
     \end{subfigure}
     
    \caption{H\'enon map: optimisation trajectories for measurements with $\sigma = 0.001$ over (a) loss and (b) error, and for measurements with $\sigma = 0.1$ over (c) loss and (d) error. Relaxation 3 (blue), Relaxation 2 (orange), and Relaxation 1 (green) for ten ground truth time series ($K = 10^5$). 
    The dashed line indicates the loss of one ground truth solution.}
    \label{figurehenon1}
\end{figure}

We compare the trajectories of the losses $\ell^{H}_{1,2,3}$ and the error $e$ for all three relaxations with $\rho = 1$ and  measurement noise with $\sigma = 0.001$ and $\sigma = 0.1$ (see Figure \ref{figurehenon1}). Both Relaxation 2 and 3 find good minima with low loss and low error in ten out of ten cases, independent of the noise level, whereas multiple trajectories of Relaxation 1 become stuck in local minima.

\subsection{State recovery of the Ikeda map}

We show a similar case with the Ikeda map
\begin{equation}
\label{ikeda1}
\begin{aligned}
   x_{1,t+1} &=& \Theta_1 + \Theta_2 (x_{1,t} \cos(\tau_t) - x_{2,t} \sin(\tau_t)), \\
   x_{2,t+1} &=& \Theta_2 (x_{1,t} \sin(\tau_t) + x_{2,t} \cos(\tau_t)),  
   \end{aligned} 
\end{equation}
and
\begin{equation}
\label{ikeda2}
\begin{aligned}
   \tau_t &=& \Theta_3 - \frac{\Theta_4}{1 + x_{1,t}^2 + x_{2,t}^2},
   \end{aligned} 
\end{equation}
where we define the measurements to be 
\begin{equation}
\label{ikeda3}
\begin{aligned}
   y_t = h(x_t) + \mu_t = ({x_{1,t}+ \tau_{t}}, {x_{2,t}+ \tau_{t}})^T + (\mu_{1,t},\mu_{2,t})^T
\end{aligned} 
\end{equation}
with $\mu_1, \mu_2 \sim \mathcal{N}(0,\sigma)$. We simulate trajectories with ground truth $\Theta_1 = 1, \Theta_2 = 0.9, \Theta_3 = 0.4,$ and $\Theta_3 = 6$.

In Relaxation 3, we reconstruct all time-varying states and parameters $x_t = (x_{1,t}, x_{2,t}, \tau_t, \theta_{1,t}, \theta_{2,t}, \theta_{3,t}, \theta_{4,t})^T$ from $y_t$, so that, with $n+m = 7$, 
\begin{equation}
\label{ikedaloss}
\begin{aligned}
    &\ell_3^{I}(\hat{\mathcal{X}}_3) = \rho \sum_{t=1}^{N}\|y_t-h(\hat{x}_{t})\|^2_2 +
    \delta \sum_{t=1}^{N-1} \|\hat{\theta}_{t+1} -\hat{\theta}_{t}\|^2_2 + \\ &\sum_{t=1}^{N-1} \| \hat{x}_{1,t+1}-(\hat{\theta}_{1,t} + \hat{\theta}_{2,t} (\hat{x}_{1,t} \cos(\hat{\tau}_t) - \hat{x}_{2,t} \sin(\hat{\tau}_t)))\|_2^2\ + \\
    &\sum_{t=1}^{N-1} \| \hat{x}_{2, t+1}-(\hat{\theta}_{2,t} (\hat{x}_{1,t} \sin(\hat{\tau}_t) + \hat{x}_{2,t} \cos(\hat{\tau}_t)))\|_2^2\ + \\
    &\sum_{t=1}^{N-1} \| \hat{\tau}_{t+1}-(\hat{\theta}_{3,t} - \frac{\hat{\theta}_{4,t}}{1 + \hat{x}_{1,t}^2 + \hat{x}_{2,t}^2})\|_2^2\ .
\end{aligned}
\end{equation}

Loss $\ell_2^{H}(\hat{\mathcal{X}}_2)$ can again be derived by replacing all time-varying $\hat{\theta}_t$ by a corresponding single estimate $\hat{\Theta}$.
In Relaxation 1, we reconstruct states $\hat{\mathcal{X}}_1$, given $\Theta_{1}, \Theta_{2}, \Theta_{3}, \Theta_{4}$ and $\hat{\tau}_t$ simply calculated from $\hat{x}_{1,t}, \hat{x}_{2,t}$ in every optimisation step according to \ref{ikeda2},
\begin{equation}
\label{ikedaloss2}
\begin{aligned}
    &\ell_1^{I}(\hat{\mathcal{X}}_1) = \rho \sum_{t=1}^{N}\|y_t-h(\hat{x}_{t})\|^2_2 + \\ 
    & \sum_{t=1}^{N-1} \| \hat{x}_{1,t+1}-({\Theta}_{1} + {\Theta}_{2} (\hat{x}_{1,t} \cos(\hat{\tau}_t) - \hat{x}_{2,t} \sin(\hat{\tau}_t)))\|_2^2\ + \\
    &\sum_{t=1}^{N-1} \| \hat{x}_{2, t+1}-({\Theta}_{2} (\hat{x}_{1,t} \sin(\hat{\tau}_t) + \hat{x}_{2,t} \cos(\hat{\tau}_t)))\|_2^2\ , 
\end{aligned}
\end{equation}
with
\begin{equation}
\begin{aligned}
   \hat{\tau}_t &=& \Theta_3 - \frac{\Theta_4}{1 + \hat{x}_{1,t}^2 + \hat{x}_{2,t}^2}.
   \end{aligned} 
\end{equation}

\begin{figure}[]
     \centering
     \begin{subfigure}{0.35\columnwidth}
    \includegraphics[width=\columnwidth]{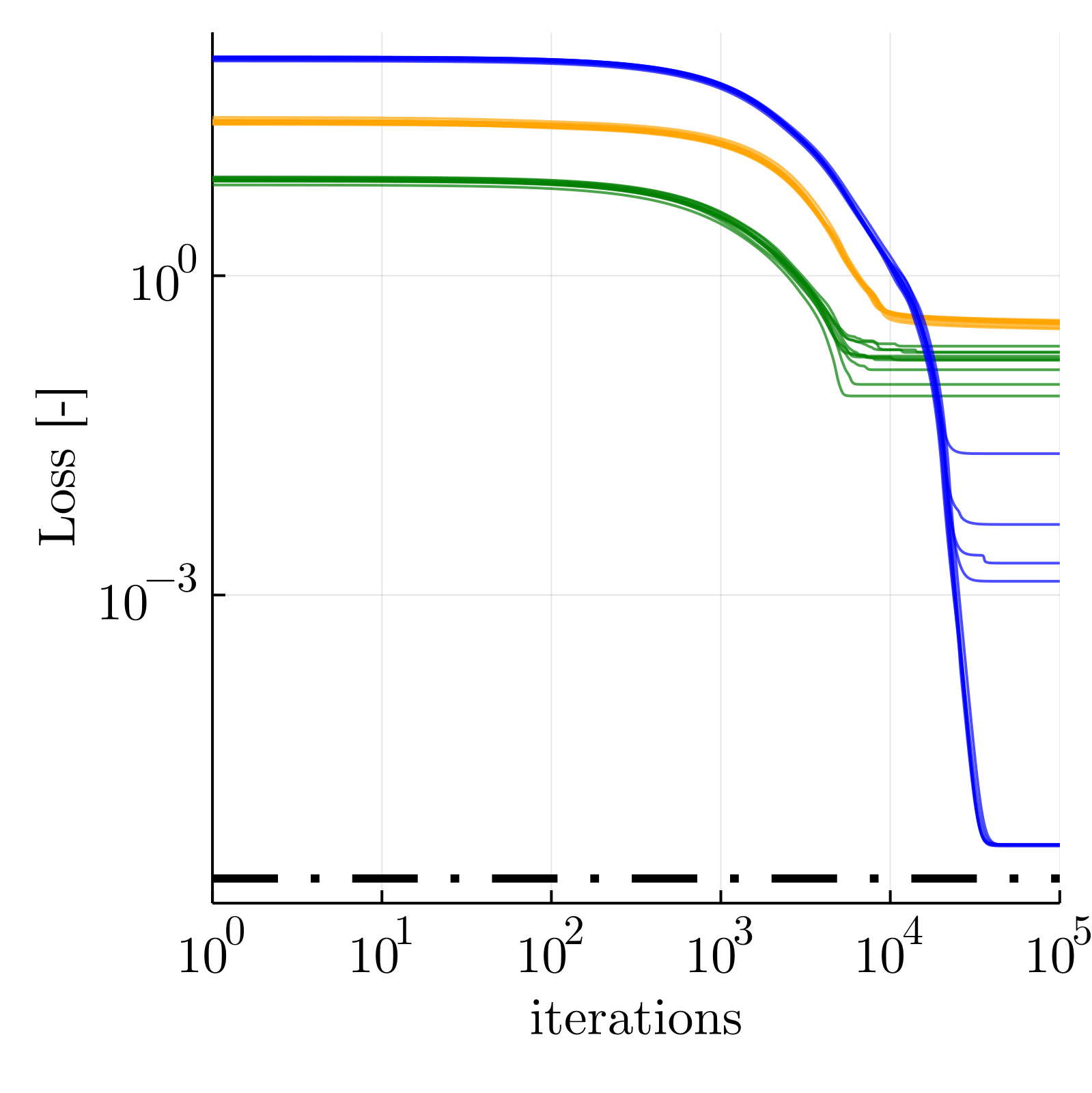}
         \caption{}
     \end{subfigure}
     \hspace{1mm}
     \begin{subfigure}{0.35\columnwidth}
    \includegraphics[width=\columnwidth]{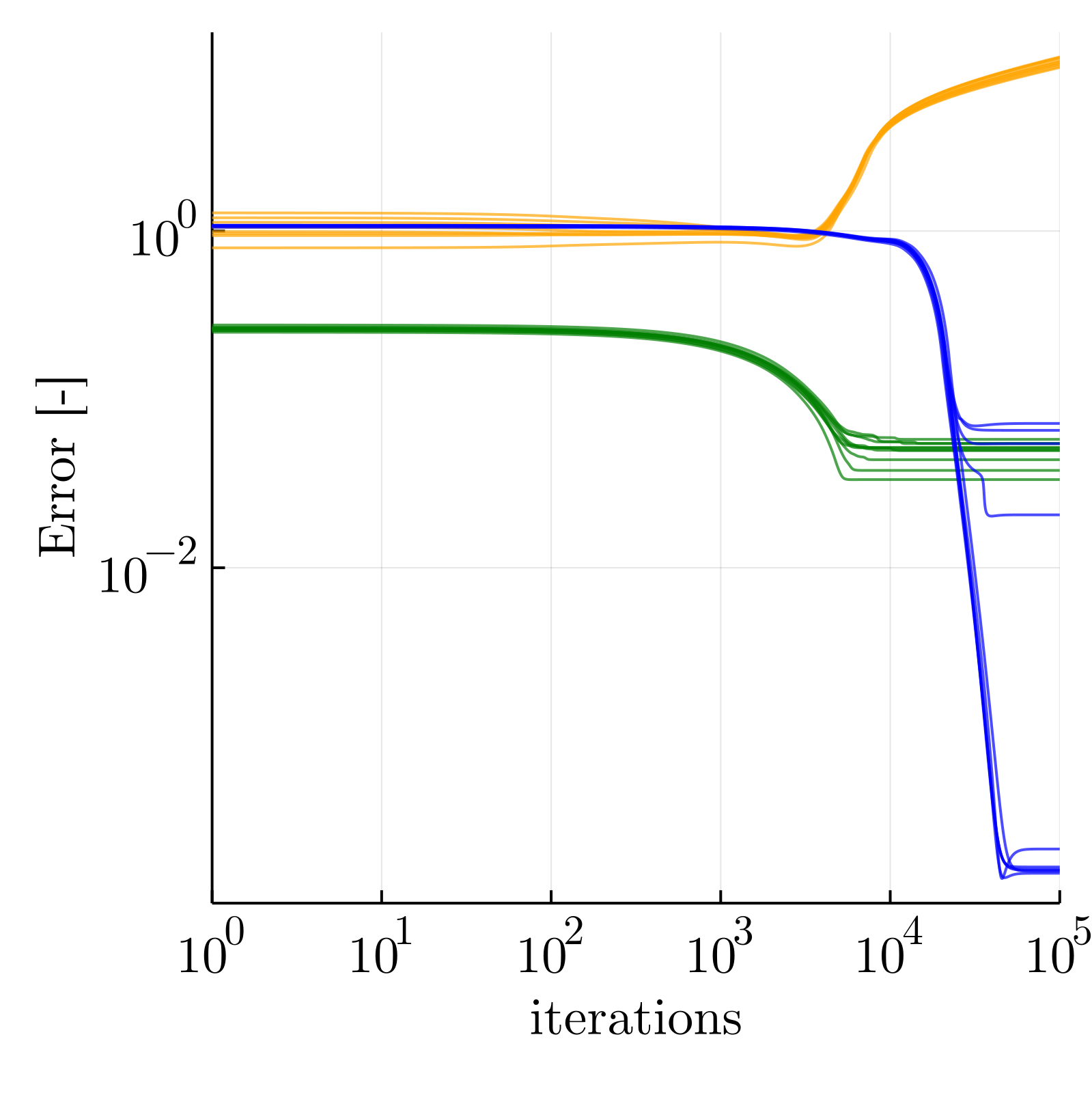}
         \caption{}
     \end{subfigure}
     \vspace{1mm}
     \begin{subfigure}{0.35\columnwidth}
    \includegraphics[width=\columnwidth]{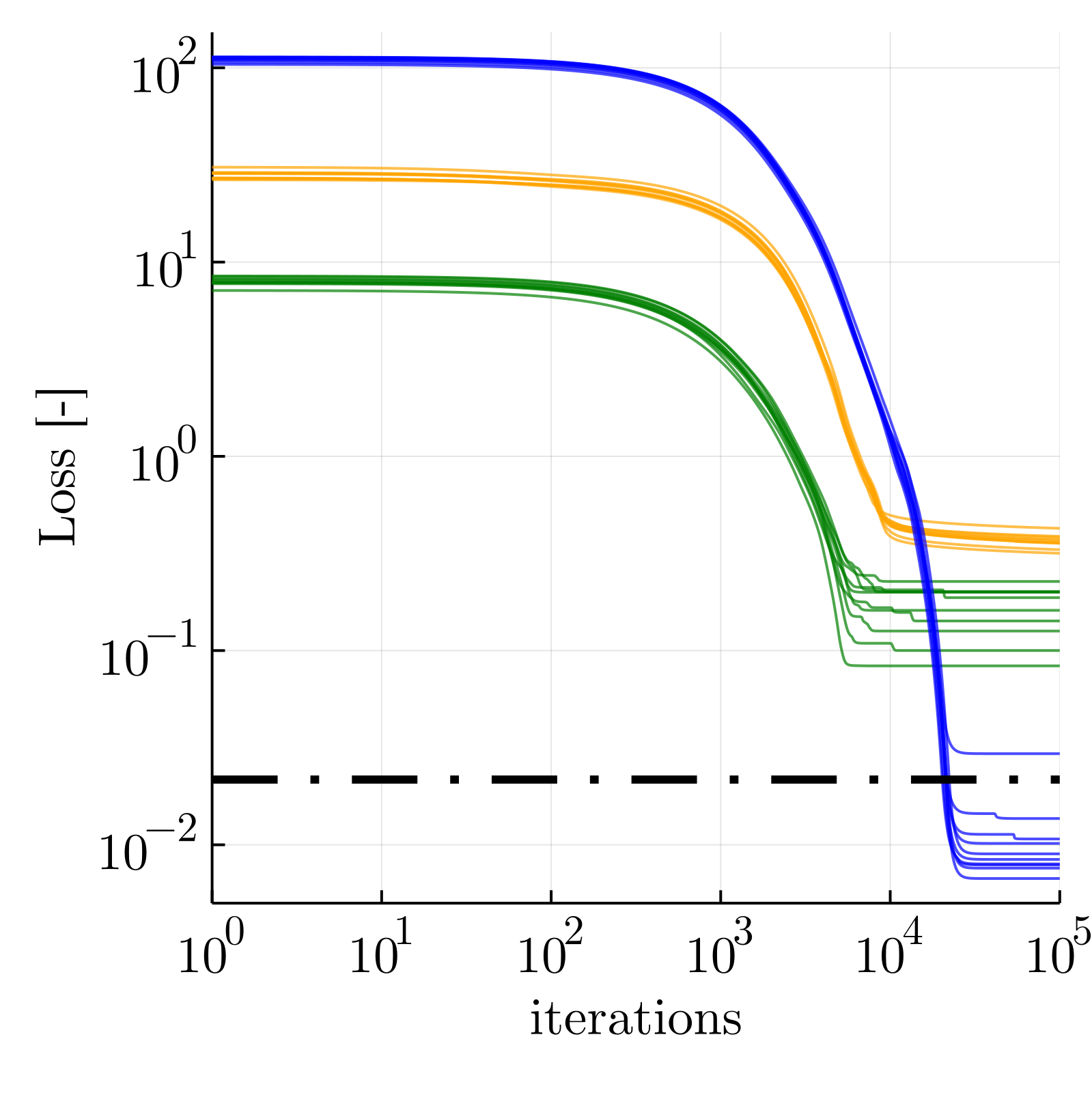}
         \caption{}
     \end{subfigure}
     \hspace{1mm}
     \begin{subfigure}{0.35\columnwidth}
    \includegraphics[width=\columnwidth]{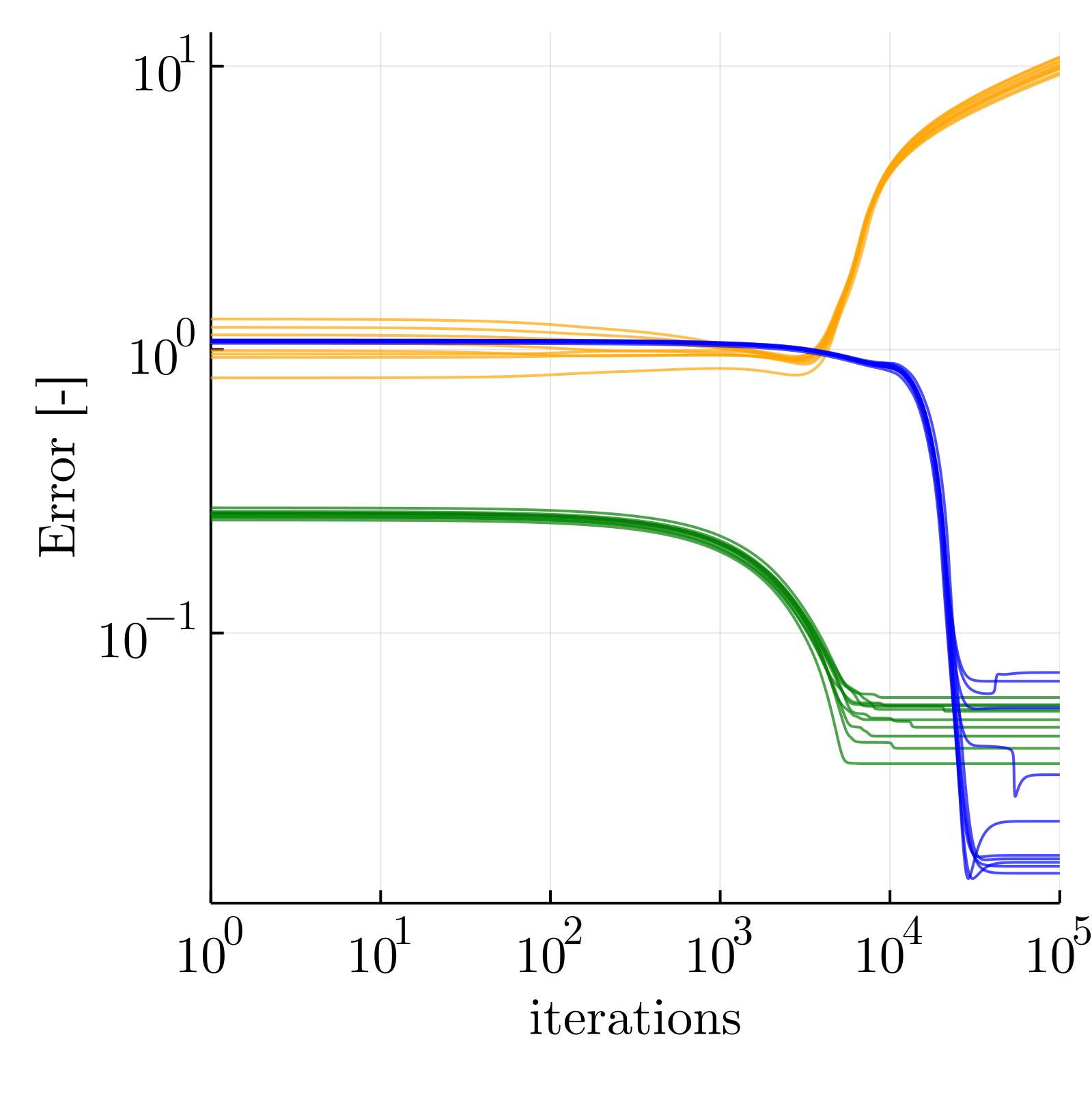}
         \caption{}
     \end{subfigure}
     
    \caption{Ikeda map: optimisation trajectories for measurements with $\sigma = 0.001$ over (a) loss and (b) error, and with $\sigma = 0.1$ over (c) loss and (d) error. Relaxation 3 (blue), Relaxation 2 (orange), and Relaxation 1 (green) each for the same ten ground truth time series ($K = 10^5$). 
    The dashed line indicates the loss of one ground truth solution.}
    \label{figureikeda1}
\end{figure}

We depict the optimisation process for $\rho = 1$ in Figure \ref{figureikeda1} for $\sigma$ equal to 0.001 and 0.1, respectively. 
Only the overparameterised loss (Relaxation 3) converges towards the ground truth in six out of 10 cases. 
Yet, we only require a reasonable number of opportunities to find solutions with low costs, as we can always explore the loss landscape from multiple starting points $\hat{\mathcal{X}}^1$. 
In contrast, all other trajectories of Relaxation 1 and 2 become stuck in local minima. 
The parameter estimates in Relaxation 2 diverge from the ground truth, leading to an increase in error over iterations.

\subsection{State recovery of the Chialvo map}

The third example is the Chialvo map,
\begin{equation}
\label{chialvox}
\begin{aligned}
   x_{1,t+1} &=&  x_{1,t}^2 \exp(x_{2,t}-x_{1,t}) + u_t, \\
   x_{2,t+1} &=& \Theta_1 x_{2,t} - \Theta_2 x_{1,t} + \Theta_3,  
   \end{aligned} 
\end{equation}
and, with $\mu_t \sim \mathcal{N}(0, \sigma)$,
\begin{equation}
\label{chialvoy}
\begin{aligned}
   y_t &=& x_{1,t} + \mu_t,
   \end{aligned} 
\end{equation}
in the chaotic range $\Theta_1 = 0.89, \Theta_2 = 0.18, \Theta_3 = 0.28$. 
We assume, apart from the measurements $y_t$, a known input $u_t = 0.025 ~ \forall t$ (for observability). 
Simulations are performed from uniformly distributed initial conditions in the oscillating regime between 1 and 2, i.e., $\kappa \sim \mathcal{U}(1,2)$.

Relaxation 3 reconstructs $x_t = (x_{1,t}, x_{2,t}, \theta_{1,t}, \theta_{2,t}, \theta_{3,t})^T$, $\forall t$ ($n+m = 5$), via the loss
\begin{equation}
\label{chialvoloss}
\begin{aligned}
    \ell_3^{C}(\hat{\mathcal{X}}_3) = &\sum_{t=1}^{N-1} \| \hat{x}_{1,t+1}-(\hat{x}_{1,t}^2 \exp(\hat{x}_{2,t}-\hat{x}_{1,t}) + u_t)\|_2^2\ + \\
    &\sum_{t=1}^{N-1} \| \hat{x}_{2, t+1}-(\hat{\theta}_{1,t} \hat{x}_{2,t} - \hat{\theta}_{2,t} \hat{x}_{1,t} + \hat{\theta}_{3,t})\|_2^2\ + \\
    &\rho \sum_{t=1}^{N}\|y_t-\hat{x}_{1,t}\|^2_2 +
    \delta \sum_{t=1}^{N-1} \|\hat{\theta}_{t+1} -\hat{\theta}_{t}\|^2_2.
\end{aligned}
\end{equation}

We then modify Relaxation 2, akin to classical system identification, so that measured states $\hat{x}_{1,t}$ are replaced with $y_t$ directly due to our knowledge of $h$. Thus,
\begin{equation}
\label{chialvoloss2}
\begin{aligned}
    \ell_2^{C}(\hat{\mathcal{X}}_2) = &\sum_{t=1}^{N-1} \| {y}_{t+1}-(y_t^2 \exp(\hat{x}_{2,t}-y_t) + u_t)\|_2^2
    + \\
    &\sum_{t=1}^{N-1} \| \hat{x}_{2,t+1}-(\hat{\Theta}_1 \hat{x}_{2,t} - \hat{\Theta}_2 y_t + \hat{\Theta}_3)\|_2^2,
\end{aligned}
\end{equation}
with estimates $\hat{\mathcal{X}}_2 = (\hat{x}_{2,1}, \cdots, \hat{x}_{2,N}, \hat{\Theta}_{1},\hat{\Theta}_{2},\hat{\Theta}_{3} )^T$. 
The task is thus to reconstruct only a state $x_{2,t} ~ \forall t$ as well as the parameters.
The explicit use of $y_t$ in the optimisation, without an additional, denoised state $\hat{x}_1$, is a common approach within the prediction error framework \cite{Schoukens2019,Brunton2016b}. 

The loss $\ell_1^{C}(\hat{\mathcal{X}}_1)$ for Relaxation 1 derives from the equation (\ref{chialvoloss2}) simply by replacing the unknown parameter estimates 
by their respective ground truth values; hence, $\hat{\mathcal{X}}_1 = (\hat{x}_{2,1}, \cdots, \hat{x}_{2,N})^T$. 

\begin{figure}[]
     \centering
     \begin{subfigure}{0.35\columnwidth}
    \includegraphics[width=\columnwidth]{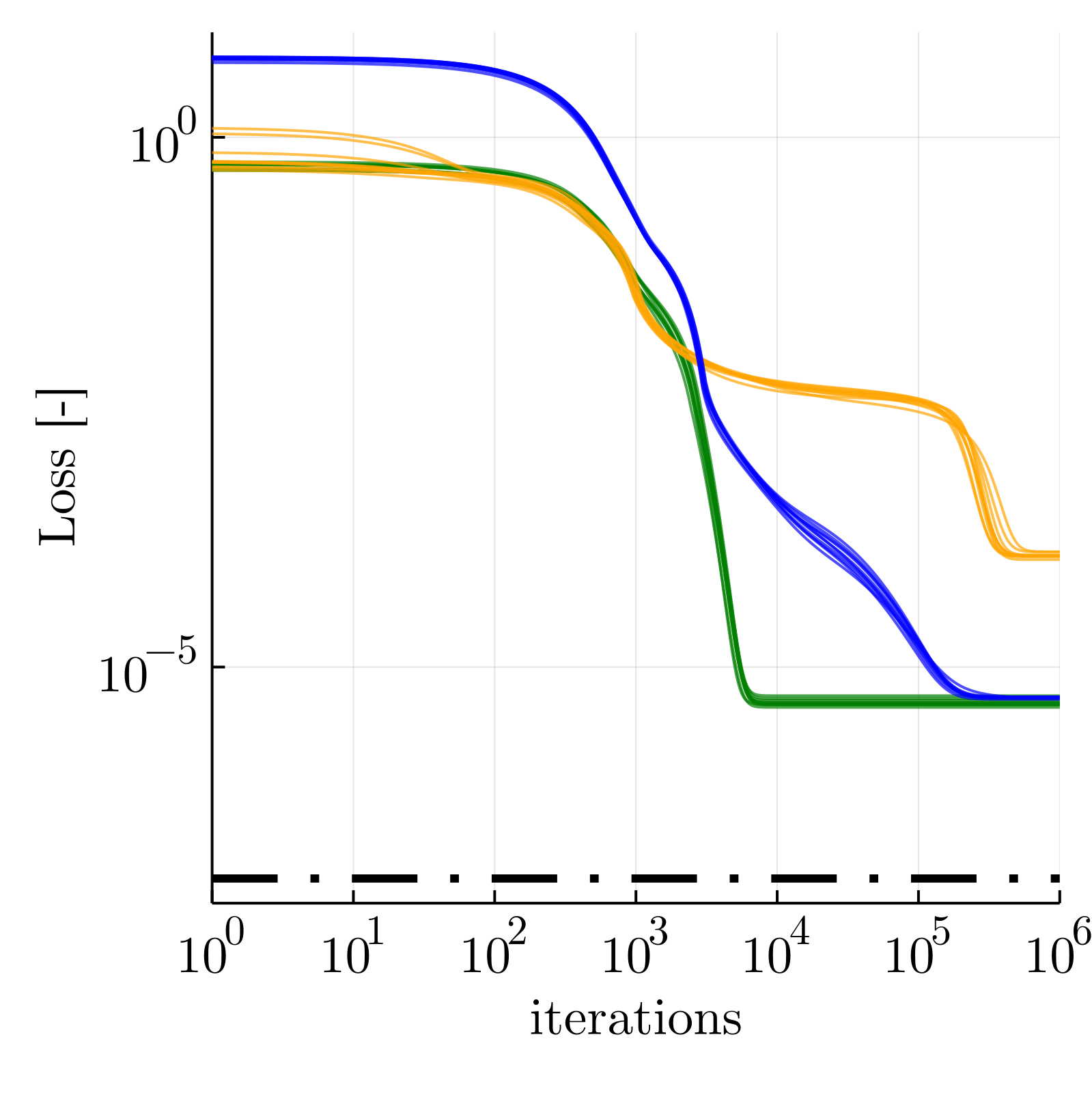}
         \caption{}
     \end{subfigure}
     \hspace{1mm}
     \begin{subfigure}{0.35\columnwidth}
    \includegraphics[width=\columnwidth]{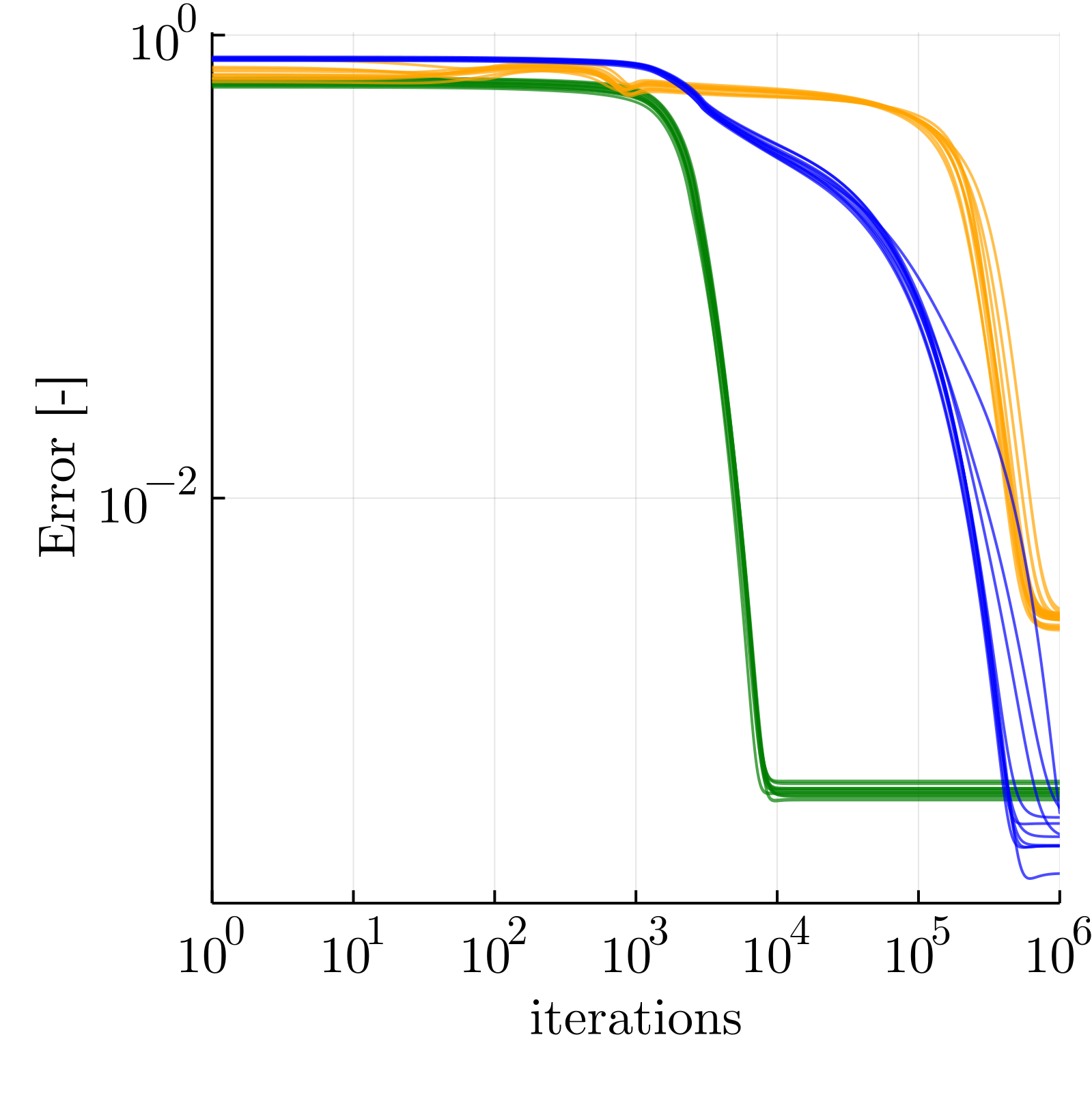}
         \caption{}
     \end{subfigure}
          \vspace{1mm}
     \begin{subfigure}{0.35\columnwidth}
    \includegraphics[width=\columnwidth]{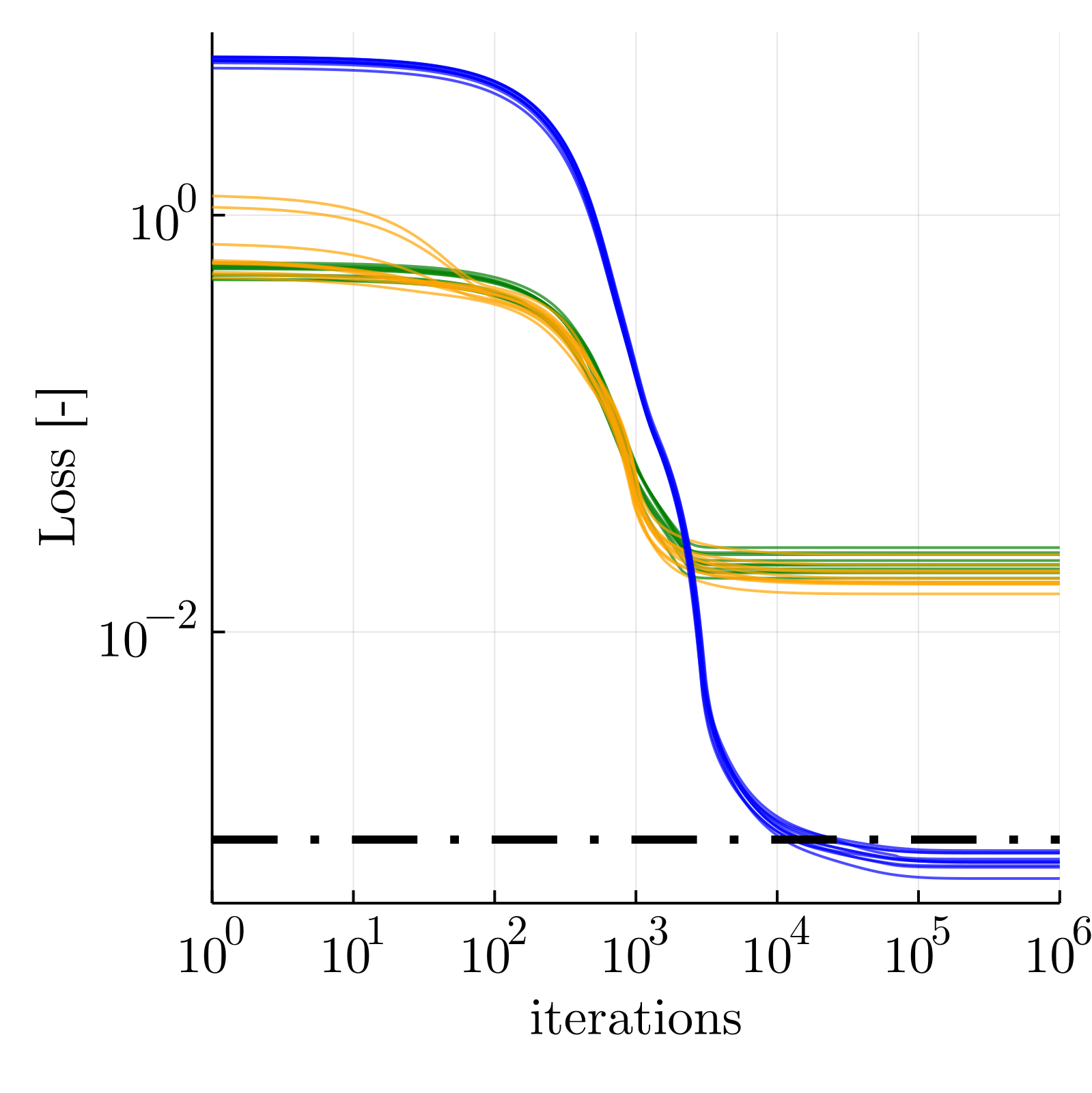}
         \caption{}
     \end{subfigure}
     \hspace{1mm}
     \begin{subfigure}{0.35\columnwidth}
    \includegraphics[width=\columnwidth]{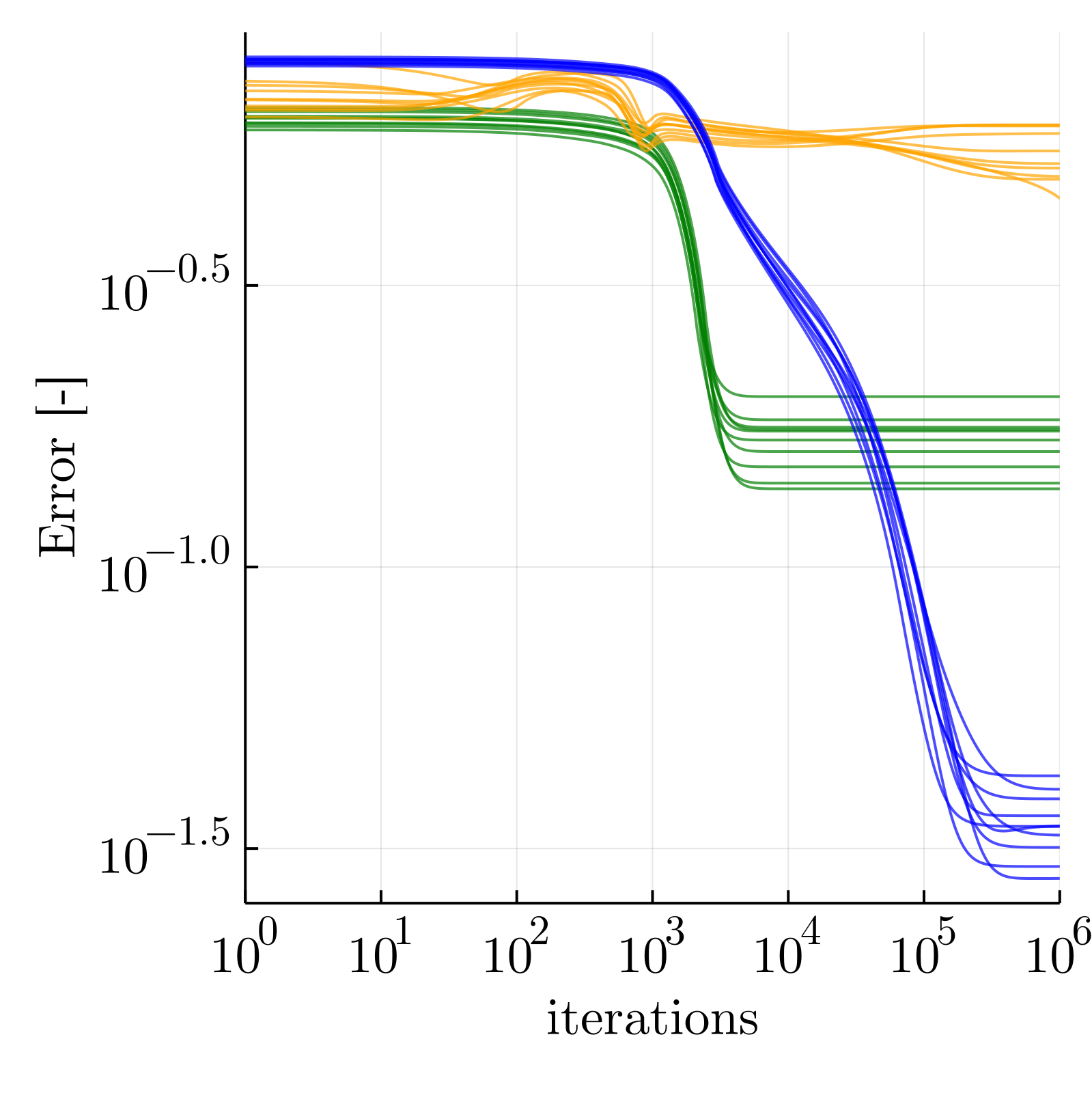}
         \caption{}
     \end{subfigure}
    \caption{Chialvo map: optimisation trajectories for measurements with $\sigma = 0.001$ over (a) loss and (b) error, and with $\sigma = 0.1$ over (c) loss and (d) error. Relaxation 3 (blue), Relaxation 2 (orange), and Relaxation 1 (green) each for the same ten ground truth time series ($K = 10^6$). 
    The dashed line indicates the loss of one ground truth solution.}
    \label{figurechialvo}
\end{figure}

Assuming prior knowledge, all states and parameters in the above equations (\ref{chialvox}) must be positive. Thus, we initialise $\hat{\mathcal{X}}^1$ with random values drawn from a uniform distribution $\mathcal{U}(0, 1)$. 
We also incorporate that the measurements are affected by noise, and thus use $\rho = 0.1$ in Relaxation 3.
See Figure \ref{figurechialvo} for results: Relaxation 3 again finds solutions closest to the ground truth, whereas Relaxation 2 converges towards local minima. Relaxation 1 works well in the case without noise, but is not as robust against higher level of measurement noise as Relaxation 3 incorporating denoising.

\begin{table}
\caption{Median computation time (in $10^{-5}$ sec.) of each $\nabla \ell(\hat{\mathcal{X}})$ over $10^6$ random $\hat{\mathcal{X}}$, each for $N=300$ as in the examples, and a scaled version with $N=1200$.}
\begin{center}
\begin{tabular}{|c||c| c | c |}
 \hline
loss & $\nabla \ell_1^H(\hat{\mathcal{X}})$ & $\nabla \ell_2^H(\hat{\mathcal{X}})$ & $\nabla \ell_3^H(\hat{\mathcal{X}})$\\
\hline
$N = 300$ & 8.5 & 8.6 & 19.5\\
\hline
$N = 1200$ & 13.0 & 14.2 & 33.9\\
\hline
loss &$\nabla \ell_1^I(\hat{\mathcal{X}})$& $\nabla \ell_2^I(\hat{\mathcal{X}})$&$\nabla \ell_3^I(\hat{\mathcal{X}})$\\
\hline
$N = 300$ & 28.9 & 27.5 & 61.1\\
\hline
$N = 1200$ & 60.0 & 78.1 & 124.6\\
\hline
loss &$\nabla \ell_1^C(\hat{\mathcal{X}})$& $\nabla \ell_2^C(\hat{\mathcal{X}})$&$\nabla \ell_3^C(\hat{\mathcal{X}})$\\
\hline
$N = 300$ & 6.5 & 6.9 & 22.1\\
\hline
$N = 1200$ & 16.7 & 18.8 & 53.8\\
\hline
\end{tabular}
\end{center}
\label{TableProj2bRes}
\end{table}

\newpage
\section{Discussion and Conclusions}

We have presented an empirical study on gradient descent-based optimisation of highly nonlinear loss functions. These arise from problems in system identification where the reconstructions of unknown system states and parameters of dynamical systems are of interest. 
So far, problems like these have been discussed theoretically \cite{Mareels2002} and, more extensively, in the linear case \cite{Haderlein1}. 
This study presents a follow-up on this previous work, through a number of nonlinear case studies. 

We have compared the optimisation of an overparameterised loss with seemingly simpler loss functions with fewer unknowns.
Surprisingly, we find that the former optimisations in higher dimensions consistently converge towards solutions with low error. 
In contrast, the cases with fewer degrees of freedom have a higher probability of becoming stuck in sub-optimal minima along the optimisation trajectories.
Our findings are robust in terms of measurement noise.
If noise is present, the overparameterised solutions can reach a loss that is potentially smaller than that of the ground truth as noise is effectively absorbed into the parameter estimates (a situation well understood in the linear case).

We want to stress that the overparameterisation in this work does not lead to interpolation of noisy measurements (`overfitting'). 
Instead, we find solutions in the unconstrained loss landscape that are very close to the true underlying dynamical system. 
A theoretical convergence assessment for such overparameterised loss functions might be possible via certain constraints on $f,h$ (e.g., controlling the Hessian tensor to achieve a well-conditioned system), enabling further insight into the relationship with machine learning. 
In contrast to such tasks that search within a class of universal approximators, we limit the space of models severely.
Therefore, this study is confined to scenarios with prior knowledge of the system (in that the desired system belongs to a given class of systems), which is in practice not to be expected.
Relaxing this condition is to be explored in future work.


\bibliographystyle{unsrt}  
\bibliography{references}

\end{document}